\documentclass[3p]{elsarticle}

\usepackage{lineno,stmaryrd}

\journal{}

\usepackage[utf8]{inputenc}
\usepackage{mathrsfs}
\usepackage{amsfonts}
\usepackage{amsthm, amsmath}
\usepackage{amssymb}
\usepackage{amscd}
\usepackage{amstext}
\usepackage[all]{xy}
\xyoption{2cell}
\UseTwocells
\usepackage{graphicx}
\usepackage{comment}
\usepackage{calrsfs}
\usepackage{framed}

\usepackage[colorinlistoftodos]{todonotes}
\usepackage[colorlinks=true, allcolors=blue]{hyperref}

\usepackage{tikz}
\usetikzlibrary{positioning,intersections,cd,matrix, arrows, backgrounds}
\usetikzlibrary{shapes,shapes.geometric,shapes.misc}
\usetikzlibrary{decorations}
\usetikzlibrary{decorations.pathreplacing}
\def\Ar#1#2#3{\ar[from={#1}, to={#2}, #3]}
\usepackage[OT2,T1]{fontenc}

\theoremstyle{plain}
\newtheorem{thm}{Theorem}[section]
\newtheorem{prp}[thm]{Proposition}
\newtheorem{lem}[thm]{Lemma}
\newtheorem{cor}[thm]{Corollary}

\newtheorem{cnj}[thm]{Conjecture}

\newtheorem*{thm-nn}{Theorem}
\newtheorem*{prp-nn}{Proposition}
\newtheorem*{lem-nn}{Lemma}
\newtheorem*{cor-nn}{Corollary}
\newtheorem*{clm-nn}{Claim}
\newtheorem*{cnj-nn}{Conjecture}
\newtheorem*{prb-nn}{Problem}

\theoremstyle{definition}
\newtheorem{dfn}[thm]{Definition}
\newtheorem{exm}[thm]{Example}

\newtheorem*{dfn-nn}{Definition}

\newtheorem{rmk}[thm]{Remark}
\newtheorem{ntn}[thm]{Notation}

\newcommand{\xyR}[1]{%
  \xydef@\xymatrixrowsep@{#1}}

\newcommand{\xyC}[1]{%
  \xydef@\xymatrixcolsep@{#1}}

\def\al{\alpha}
\def\be{\beta}

\def\de{\delta}

\def\ro{\rho}

\def\ta{\tau}

\def\Ga{\Gamma}

\def\La{\Lambda}

\def\Ker{\operatorname{Ker}}

\def\Im{\operatorname{Im}}

\def\Hom{\operatorname{Hom}}

\def\rad{\operatorname{rad}}
\def\Aut{\operatorname{Aut}}
\def\End{\operatorname{End}}
\def\Ext{\operatorname{Ext}}

\def\mod{\operatorname{mod}}

\def\supp{\operatorname{supp}}

\def\add{\operatorname{add}}

\def\calC{{\mathcal C}}

\def\calE{{\mathcal E}}

\def\calI{{\mathcal I}}

\def\calL{{\mathcal L}}

\def\calP{{\mathcal P}}

\def\calS{{\mathcal S}}

\def\calX{{\mathcal X}}

\def\bbN{{\mathbb N}}
\def\bbZ{{\mathbb Z}}

\def\bbI{{\mathbb I}}
\def\bbJ{{\mathbb J}}
\def\op{^{\mathrm{op}}} 

\def\inv{^{-1}}

\def\iso{\cong}
\def\ds{\oplus}

\def\ovl{\overline}
\def\Ds{\bigoplus}
\def\DS{\bigoplus\limits}

\def\dsm#1,#2..#3{\bigoplus_{{#1}={#2}}^{#3}}
\def\sm#1,#2..#3{\sum_{{#1}={#2}}^{#3}}

\def\id{1\kern-.25em{\text{{\rm l}}}} 

\def\isoto{\ \raise.8ex\hbox{$^{\sim}$}\kern-.7em\hbox{$\to$}\ }
\def\Cdot{\centerdot}
\def\down{_{\Cdot}}

\def\ya#1{\xrightarrow{#1}}
\def\blank{\operatorname{-}}

\def\bg{%
  \family{cmr}\size{20}{12pt}\selectfont}

\def\bigzerou{%
  \smash{\lower1.7ex\hbox{\bg 0}}}

\def\repr[#1;#2;#3;#4;#5]{
  \left(
    \begin{matrix}#1\\#2\end{matrix}
    #3
    \begin{matrix}#4\\#5\end{matrix}
  \right)}
\def\spmat#1{\left(\begin{smallmatrix} #1 \end{smallmatrix}\right)}

\usepackage[e]{esvect}

\newcommand{\Gf}[1]{\vv*{G}{#1}}

\def\k{\Bbbk}

\newcommand{\bi}[3]{{}_{#2}{#1}_{#3}}
\newcommand{\hyph}{\text{-}}

\newcommand{\cov}{\operatorname{Cov}}

\newcommand{\vect}{\mathrm{vect}}

\newcommand{\cone}{\operatorname{cone}}
\newcommand{\rdim}{\operatorname{res-dim}}
\newcommand{\rgldim}{\operatorname{res-gldim}}
\newcommand{\gldim}{\operatorname{gldim}}
\newcommand{\pd}{\operatorname{pd}}
\newcommand{\intdim}{\operatorname{int-res-dim}}

\newcommand\comcat{\!\downarrow\!}
\newcommand\intgldim{\operatorname{int-res-gldim}}

\DeclareMathOperator{\Seg}{Seg}
\newcommand\intval[1]{\llbracket #1 \rrbracket}

\begin{document}

\begin{frontmatter}
  \title{Approximation by interval-decomposables and interval resolutions of persistence modules%
    \tnoteref{titlededication}\tnoteref{mytitlenote}}

  \tnotetext[titlededication]{
    Dedicated to the memory of Daniel Simson and Andrzej Skowro{\'n}ski
  }

  \tnotetext[mytitlenote]{
    H.A. is supported by JSPS Grant-in-Aid for Scientific Research (C) 18K03207, and
    JSPS Grant-in-Aid for Transformative Research Areas (A) (22A201).
    E.G.E. is supported by JSPS Grant-in-Aid for Transformative Research Areas (A) (22H05105).
    K.N. is supported by JSPS Grant-in-Aid for Transformative Research Areas (A) (20H05884).
    M.Y. is supported by JSPS Grant-in-Aid for Scientific Research (C) (20K03760).
    H.A. and M.Y.\ are partially supported by Osaka Central Advanced Mathematical Institute (JPMXP0619217849).}

  \author[shizu,kuias,ocami]{Hideto Asashiba}
  \ead{asashiba.hideto@shizuoka.ac.jp}

  \author[kobe]{Emerson G. Escolar\corref{mycorrespondingauthor}}
  \cortext[mycorrespondingauthor]{Corresponding author}
  \ead{e.g.escolar@people.kobe-u.ac.jp}

  \author[okayama]{Ken Nakashima}
  \ead{ken.nakashima@okayama-u.ac.jp}

  \author[ocami]{Michio Yoshiwaki}
  \ead{yosiwaki@sci.osaka-cu.ac.jp, michio.yoshiwaki@omu.ac.jp}

  \address[shizu]{Faculty of Science, Shizuoka University, Japan}
  \address[kuias]{Institute for Advanced Study, Kyoto University, Japan}
  \address[ocami]{Osaka Central Advanced Mathematical Institute, Osaka, Japan}
  \address[kobe]{Graduate School of Human Development and Environment, Kobe University, Japan}
  \address[okayama]{Center for AI, Mathematical and Data Sciences, Okayama University, Japan}

  \begin{abstract}
    In topological data analysis, two-parameter persistence can be studied using the representation theory of the $2$d commutative grid, the tensor product of two Dynkin quivers of type A. In a previous work, we defined interval approximations using restrictions to essential vertices of intervals together with Mobius inversion.
    In this work, we consider homological approximations using interval resolutions, and show that the interval resolution global dimension is finite for finite posets and that it is equal to the maximum of the interval dimensions of the Auslander-Reiten translates of the interval representations.
    In fact, for the latter equality, we obtained a general formula in the setting of finite-dimensional algebras and resolutions relative to a generator-cogenerator.
    Furthermore, in the commutative ladder case, by a suitable modification of our interval approximation, we provide a formula linking the two conceptions of approximation.
  \end{abstract}

\end{frontmatter}




\section{Introduction}

The field of topological data analysis is a rapidly growing field
applying the ideas of algebraic topology for data analysis.
In particular, one of its main tools,
persistent homology \cite{edelsbrunner2002topological,edelsbrunner2008persistent}, is able to express the
``birth'' and ``death'' of topological features in data with respect to some parameter.
The so-called interval modules play a central role in persistent homology,
as they are used to express the ``birth'' and ``death'' of topological features.
Algebraically, this follows from the fact that every pointwise finite-dimensional module of a
totally ordered poset decomposes as a direct sum of interval modules \cite{gabriel1992algebra,botnan2020decomposition,crawley2015decomposition}. For details concerning one-parameter persistent homology, we refer the reader to
\cite{edelsbrunner2002topological,edelsbrunner2008persistent}

However, when dealing with multiple parameters, as in multiparameter persistent homology \cite{carlsson2009theory},
there are difficulties with adopting exactly the same techniques.
In particular, there are many indecomposable modules that are not interval modules (see \cite{buchet_et_al:socg,buchet2020every} for concrete examples). In terms of representation theory, this means that for a large enough underlying grid of parameters (an $n$-dimensional commutative grid for $n$ parameters), the corresponding algebra is of wild representation type (\cite[Theorem~1.3]{bauer2020cotorsion}, \cite[Theorem~2.5]{leszczynski1994representation}, \cite[Theorem~5]{leszczynski2000tame}).

Thus there are many attempts to construct new invariants to capture persistent topological information even in the multiparameter setting.
One idea is to use the intervals, as they have served very well in the one-parameter case.
In particular, one recent direction is to study generalized persistence diagrams
\cite{kim2018generalized,asashiba2019approximation,botnan2021signed,dey2021computing}
which associates ``multiplicities'' to intervals or subclasses of intervals,
much like how the persistence diagram $d_M$ of a persistence module $M$ associates to each interval $I$
its multiplicity $d_M(I)$ as a direct summand in $M$.

Recently, the works \cite{botnan2021signed,blanchette2021homological} introduced the explicit use of
relative homological algebra in persistence theory. With this, a new perspective on
several existing invariants in multiparameter persistence was provided.
In particular, one of the contributions of \cite{botnan2021signed}
is the study of rank decompositions of the rank invariant,
and their relationships with the generalized persistence diagrams and with projective resolutions relative to the so-called rank-exact structure.
The work \cite{blanchette2021homological} provides a systematic use of relative homological
algebra in the study of invariants via their framework of ``homological invariants''.

In this work, we 
study the relative homological algebra with respect to the interval modules.
In particular, we show that for an arbitrary finite poset, the interval resolution global dimension of its incidence algebra is finite (Proposition~\ref{prp:intgldim-fin}).
This gives an affirmative answer to  Question~6.5 (whether or not the set of all interval modules $\bbI$ gives a homological invariant) of \cite{blanchette2021homological}.

We then provide a formula to compute the interval resolution global dimension as the maximum of a finite number of terms (Proposition~\ref{prp:intdim-ta}), which we use in computational experiments (in fact, the formula is proved in a more general setting; see Proposition~\ref{prp:I-res-gldim}). In particular, we obtain some conjectures about the value of the interval resolution global dimension of the $2$D commutative grids (Conjectures~\ref{cnj:ladder}~and~\ref{cnj:stabilize}).

We previously proposed a notion of compressed multiplicities relative to
some set of essential vertices of intervals \cite{asashiba2019approximation}.
These compressed multiplicities are related to (and can be seen as variants of)
the generalized rank invariant proposed in \cite{kim2018generalized}
(see \cite{asashiba2019approximation} for a careful comparsion).
Applying a M\"obius inversion, one respectively gets a generalized persistence diagram.
However, \cite{blanchette2021homological} showed that the generalized rank invariant
(and thus the corresponding generalized persistence diagram)
is not a dim-hom invariant relative to $\bbI$, and conjectures that it is not a homological invariant.

Here, in the case of the commutative ladder (an $m \times 2$ commutative grid), we
provide a modification (Definition~\ref{def:compressedmult}) of our previous notion of compressed multiplicity.
This gives a new invariant that is intimately related to an alternating sum of the terms appearing in the interval resolution (Theorem~\ref{thm:comp-mult-resol}) and
whose M\"obius inversion is a homological invariant relative
to $\bbI$ (Corollary~\ref{cor:homologicalinvariant}).
That is, this gives an alternative definition for the homological invariant relative to $\bbI$.

We note that this paper is heavily grounded in the representation theory of algebras
and uses its methods freely.
For example,
Proposition~\ref{prp:intgldim-fin} is proved by using a Theorem \cite[Theorem~in~\S 5]{ringel2010iyama}
on quasi-hereditary algebras due to Ringel (see also Iyama's \cite[Lemma~2.2]{iyama2003finiteness}).
The proofs of Proposition~\ref{prp:intdim-ta} and Theorem~\ref{thm:comp-mult-resol}
essentially use the theory of almost split sequences due to Auslander and Reiten
(see e.g.~\cite{auslander1997representation}), and
one of the key points for the latter is in the fact that the left adjoint functor
to the restriction functor (compression functor) sends all modules (involved in an almost split sequence) to interval decomposable modules.
Thus, we assume background knowledge in the representation theory of algebras, for which
we refer the reader to \cite{assem2006elements,auslander1997representation}.


\section{Background}

Throughout this paper $\k$ is a field, and $\calP$ always denotes a {\em finite} poset.
We set $D:= \Hom_\k(\hyph, \k)$ to be the usual self-duality of the category $\mod A$ of finite-dimensional $A$-modules for a finite-dimensional $\k$-algebra $A$.
We use the following definition concerning subsets of a finite poset $\calP$.
\begin{dfn}
  \label{dfn:posetsubsets}
  Let $\calP$ be a finite poset.
  \begin{enumerate}
  \item For $p, q \in \calP$, the \emph{segment} from $p$ to $q$ is
    $[p,q] := \{x \in \calP \mid p \leq x \leq q\}$ \label{dfn:segment}
  \item A subset $S \subseteq \calP$ is said to be \emph{connected} if it is connected as a subgraph of the Hasse diagram of $\calP$
  \item A subset $S \subseteq \calP$ is said to be \emph{convex} if for any
  $p,q \in {S}
  $, $[p,q] \subseteq S$.
  \item A subset $S \subseteq \calP$ is said to be an \emph{interval} if it is connected and convex.
  \item For $p\in \calP$, a subset $\cov(p)$ is defined to be the set of the elements $q$ covering $p$, that is, $q$ satisfies the conditions $p\le q$ and $[p,q]=\{p,q\}$.  
  \item An element $p\in \calP$ is said to be an {\em upper bound} of $S$ if $s\leq p$ for each $s\in S$. 
  The set of upper bounds of $S$ is denoted by $U(S)$.
  \item An element $p\in U(S)$ is said to be the {\em join} of $S$ if $p\leq u$ for each $u\in U(S)$. 
  Note that the join of $S$ is unique if it exists, and is denoted by $\bigvee S$.
    \label{dfn:interval}
  \end{enumerate}
\end{dfn}
We use the notation $\bbI(\calP)$ for the set of all intervals in $\calP$.
Where the poset is clear, we simplify this to just $\bbI$.
Note that in the literature, it is standard to call the subsets defined in Definition~\ref{dfn:posetsubsets}(\ref{dfn:segment}) as ``intervals''. However, following the persistence literature, we use the word ``interval'' for Definition~\ref{dfn:posetsubsets}(\ref{dfn:interval}). Note that however, to avoid this potential source of confusion, \cite{blanchette2021homological} calls them ``spreads''. We do not adopt this terminology.

Let $\Seg(\calP)$ be the set of segments in $\calP$.
The \emph{incidence algebra} of {$\calP$} over $\k$ is the $\k$-vector space of functions from $\Seg({\calP})$ to $\k$,
with ``pointwise'' $+$ operation and $\k$-multiplication, together with the following convolution $*$ as the multiplication operation. For $f,g:\Seg({\calP}) \rightarrow \k$, define $f*g:\Seg({\calP}) \rightarrow \k$ by
 \[
   (f*g)([x,y]) := \sum\limits_{x\leq z \leq y} g([x,z])f([z,y]).
   \footnote{Note that this definition is opposite to the usual definition of the incidence algebra \cite{rota1964foundations}. We adopt this multiplication to make it match up to the usual interpretation of composition of functions from right to left.}
 \]
We denote by $\k\calP$ the incidence algebra of $\calP$ over $\k$.

The incidence algebra can be identified with the path algebra of a bound quiver $(Q,R)$, where $Q$ is the Hasse diagram of $\calP$ regarded as a quiver\footnote{For any $x,y \in \calP$,
  the quiver has a (unique) arrow $x \to y$ iff there is an edge between $x$ and $y$ in $H$ with $x \le y$.
In the path algebra, we compose arrows from the left. That is, for $\alpha: x\rightarrow y$ and $\beta: y\rightarrow z$, $\beta \alpha$ is the path from $x$ to $z$ going through $\alpha$ and then $\beta$, consistent with the convolution in the incidence algebra.}
and $R$ is the two-sided ideal of the path algebra of $Q$ generated by all the commutativity relations.

Let $\mod \k\calP$ be the category of finitely generated left $\k\calP$ modules.
We also consider $\calP$ as a category with a unique morphism $p \rightarrow q$ whenever $p \leq q$.
A (finite-dimensional) persistence module over $\calP$ is a functor $M : \calP \rightarrow \vect_\k$,
where $\vect_\k$ is the category of finite-dimensional $\k$-vector spaces. Note that $\mod \k\calP$ can be identified with the category of persistence modules over $\calP$,
which, of course, can be identified with the category of (finite-dimensional) representations of the
bound quiver $(Q,R)$.

Of great importance in persistence are the interval persistence modules.
\begin{dfn}[Interval persistence modules]
  \label{def:interval_rep}
  Let $\calP$ be a poset.
  A persistence module $M$ over $\calP$
  is said to be an \emph{interval persistence module} (\emph{interval module}, for short) if
  \begin{itemize}
  \item $\dim M(x) \leq 1$ for each $x\in \calP$,
  \item its ``support'' $\supp(M) :=\{x \in \calP \mid M(x) \neq 0\}$ is an interval of $\calP$, and
  \item for all $x \leq y$ with $x,y \in \supp(M)$, $M(x\rightarrow y)$ is an identity map.
  \end{itemize}
\end{dfn}
For each $I \in \bbI(\calP)$, the interval module with
support $I$ is uniquely determined up to isomorphism. We denote that interval module by $V_I$.

\begin{dfn}[$2$D commutative grids and commutative ladders]
  Let $m, n \in \bbN = \{1,2,3,\hdots\}$.
  \begin{itemize}
  \item The \emph{commutative $2$D grid of size} $m \times n$ is the bound quiver defined by
    the following quiver $Q$ together with all possible commutativity relations:
    \begin{equation}
      \label{eq:commgrid}
      \begin{tikzcd}
        (1,n) & (2,n) & \cdots & (m,n)\\
        \vdots & \vdots &\vdots& \vdots\\
        (1,2) & (2,2) & \cdots & (m,2)\\
        (1,1) & (2,1) & \cdots & (m,1)
        \Ar{1-1}{1-2}{"{(\al_1,n)}"}
        \Ar{1-2}{1-3}{"{(\al_2,n)}"}
        \Ar{1-3}{1-4}{"{(\al_{m-1},n)}"}
        \Ar{2-1}{2-2}{"\vdots", phantom}
        \Ar{3-1}{3-2}{"{(\al_1,2)}"}
        \Ar{3-2}{3-3}{"{(\al_2,2)}"}
        \Ar{3-3}{3-4}{"{(\al_{m-1},2)}"}
        \Ar{4-1}{4-2}{"{(\al_1,1)}"'}
        \Ar{4-2}{4-3}{"{(\al_2,1)}"'}
        \Ar{4-3}{4-4}{"{(\al_{m-1},1)}"'}
        \Ar{4-1}{3-1}{"{(1,\be_1)}"}
        \Ar{3-1}{2-1}{"{(1,\be_2)}"}
        \Ar{2-1}{1-1}{"{(1,\be_{n-1})}"}
        \Ar{4-2}{3-2}{"{(2,\be_1)}"'}
        \Ar{3-2}{2-2}{"{(2,\be_2)}"'}
        \Ar{2-2}{1-2}{"{(2,\be_{n-1})}"'}
        \Ar{4-3}{3-3}{"\cdots", phantom}
        \Ar{3-3}{2-3}{"\cdots", phantom}
        \Ar{2-3}{1-3}{"\cdots", phantom}
        \Ar{4-4}{3-4}{"{(m,\be_1)}"'}
        \Ar{3-4}{2-4}{"{(m,\be_2)}"'}
        \Ar{2-4}{1-4}{"{(m,\be_{n-1})}"'}
      \end{tikzcd}
    \end{equation}
    We denote this bound quiver by $\Gf{m,n} = (Q,\rho)$.
  \item From the discussion above, pfd representations of the bound quiver $\Gf{m,n}$
    can be identified with the pfd persistence modules over the poset
    given by $\{1,2,\hdots, m\} \times \{1,2,\hdots, n\}$ with partial order defined by
    $(i,j) \leq (k,\ell)$ if and only if
    $i \leq k$ and $j \leq \ell$. By abuse of notation, we also denote this poset by $\Gf{m,n}$.

  \item When $n = 2$ (or symmetrically $m=2$), $\Gf{m,n}$ is called the \emph{commutative ladder} of length $m$ (or $n$),
    which was studied in the context of persistence in \cite{escolar2016persistence}.
  \end{itemize}
\end{dfn}

It is known that each interval $I$ in $\bbI(\Gf{m,n})$ has a ``staircase'' form
(see the discussion in Section~4.1 of \cite{asashiba2022interval}).
That is, each interval $I$ of $\vec{G}_{m,n}$ is 
a full subposet induced by a set of the form
\[
  I_0 = \left\{ (j,i) \mid i \in \{s,s+1,\hdots, t\}, j \in \{b_i, b_i+1, \hdots, d_i\} \right\}
\]
for some $1\leq s\leq t \leq n$ and some $1 \leq b_i \leq d_i \leq m$ for each $s\leq i\leq t$ such that
\[
  b_{i+1}\leq b_{i}\leq d_{i+1}\leq d_{i}
\]
for all $i\in \{s,\dots,t-1\}$.
We follow the notation of \cite{asashiba2022interval} in writing
\[
  I_0  = \bigsqcup_{i=s}^t [b_i,d_i]_i
\]
to denote the set of points $I_0$ of the interval $I$ above.
In this notation, each $[b_i,d_i]_i$ is the ``slice'' of the staircase at height $i$.
We identify this set of points $I_0$ with the interval $I$.
\begin{exm}
We display posets using their Hasse diagrams.
Below is an interval $I$ (filled-in points and arrows), denoted $[5,6]_1 \sqcup [3,5]_2 \sqcup [3,4]_3$, of $\vec{G}_{6,4}$.
The corresponding interval module $V_I$ is given to its right.
{\tikzcdset{graphstyle/.append style={row sep=1.5em, column sep=1.5em, nodes={inner sep=0.5pt}, arrows={-stealth}}}
\begin{equation} \label{eq:st46}
  \newcommand{\bb}{\bullet}
  \newcommand{\rard}{\rar[dashed,gray]}
  \newcommand{\uard}{\uar[dashed,gray]}
  I:
  \begin{tikzcd}[graphstyle,every matrix/.append style={name=m}]
    \circ \rard & \circ \rard & \circ \rard & \circ \rard & \circ \rard & \circ \\
    \circ \rard\uard & \circ \rard\uard & \bb \rar\uard & \bb \rard\uard & \circ \rard\uard & \circ \uard\\
    \circ \rard\uard & \circ \rard\uard & \bb \rar\uar & \bb \rar\uar & \bb \rard\uard & \circ \uard\\
    \circ \rard\uard & \circ \rard\uard & \circ \rard\uard & \circ \rard\uard & \bb \rar\uar & \bb \uard
  \end{tikzcd}
  \;\;\;\;
  V_I:
  \begin{tikzcd}[graphstyle]
    0 \rar & 0 \rar & 0 \rar & 0 \rar & 0 \rar & 0 \\
    0 \rar\uar & 0 \rar\uar & \k \rar{1}\uar & \k \rar\uar & 0 \rar\uar & 0 \uar\\
    0 \rar\uar & 0 \rar\uar & \k \rar{1}\uar{1} & \k \rar{1}\uar{1} & \k \rar\uar & 0 \uar\\
    0 \rar\uar & 0 \rar\uar & 0 \rar\uar & 0 \rar\uar & \k \rar{1}\uar{1} & \k \uar
  \end{tikzcd}
\end{equation}
}
\end{exm}

We use the following notation.
\begin{ntn}
  For a finite-dimensional algebra $A$ over $\k$,
  fix a complete set $\calL$ of representatives of
  isoclasses of indecomposable modules in $\mod A$.
  For $X$ in $\mod A$,
  by the Krull--Schmidt theorem,
  there exists a unique function $d_X \colon \calL \to \bbZ_{\ge 0}$
  such that $X \iso \Ds_{L \in \calL}L^{(d_X(L))}$.
\end{ntn}


\section{Right approximations and resolutions}

Throughout this section $A$ is a finite-dimensional algebra over $\k$.
We here prepare necessary facts for right minimal approximations and resolutions with respect to
the subcategory $\calI = \add G$ for a generator $G$ in the category $\mod A$
of finitely generated left $A$-modules.
Note that all of the material here is not new,
as it is standard in relative homological algebra
(see e.g., \cite{erdmann2004radical}, \cite{dugas2007representation}, or \cite{ blanchette2021homological}).

We start with the following definition where $\calI$ is not necessarily given by $\add G$ for some generator $G$.
\begin{dfn}
  \label{dfn:approx}
  Let $M$ be in $\mod A$, and $\calI$ an additive subcategory of $\mod A$.
  \begin{enumerate}
  \item
    A right $\calI$-approximation of $M$ is a morphism
    $f \in \Hom_A(X, M)$ with $X \in \calI$ satisfying the condition that
    for any $g \in \Hom_A(Y, M)$ with $Y \in \calI$,
    there exists some $h \in \Hom_A(Y,X)$ such that $g = fh$.
    This condition is equivalent to saying that
    $f$ induces an epimorphism
    \[
      \Hom_A(\blank,f) \colon \Hom_A(\blank, X)|_{\calI} \to \Hom_A(\blank, M)|_{\calI}.
    \]
  \item
    An {\em $\calI$-resolution} of $M$ is a sequence
    \[
      \cdots \to X_2 \ya{f_2} X_1 \ya{f_1} X_0 \ya{f_0} M
    \]
    such that $f_0$ is a right $\calI$-approximation of $M$,
    and for each $i\ge 1$,
    $f_i$ is a right $\calI$-approximation of $\Ker f_{i-1}$.
    Therefore, this sequence yields an exact
 sequence
    \begin{equation}\label{eq:resol-ex}
      \begin{aligned}
        \cdots \to \Hom_A(\blank,X_2)|_{\calI} &\ya{\Hom_A(\blank,f_2)}
        \Hom_A(\blank,X_1)|_{\calI}\\
        &\ya{\Hom_A(\blank,f_1)}
        \Hom_A(\blank,X_0)|_{\calI} \ya{\Hom_A(\blank,f_0)}
        \Hom_A(\blank,M)|_{\calI} \to 0.
      \end{aligned}
    \end{equation}
  \end{enumerate}
  Dually we define a {\em left $\calI$-approximation} of $M$,
  and an {\em $\calI$-coresolution} of $M$.
\end{dfn}

\begin{rmk}\label{rmk:generator-ex}
  Let $M$ and $\calI$ be as in the definition above.
  Assume that $\calI$ contains $A$, or equivalently that
  $\calI$ contains all finitely generated projective $A$-modules
  (equivalently, $\calI$ is a generator).
  Then any right $\calI$-approximation $f \colon X \to M$ of $M$ is actually an epimorphism.
  Indeed, let $g \colon P \to M$ be an epimorphism from a finitely generated projective module $P$
  (e.g., a projective cover of $M$).
  Then since $P \in \calI$, we have $g = fh$ for some $h\colon P \to X$,
  which shows that $f$ is an epimorphism.
  Therefore, any $\calI$-resolution of $M$ in (2)
  above turns out to be an exact sequence
  \[
    \cdots \to X_2 \ya{f_2} X_1 \ya{f_1} X_0 \ya{f_0} M \to 0.
  \]

The dual remark holds for left $\calI$-approximations and $\calI$-coresolutions. Namely, if $\calI$ contains all indecomposable injective $A$-modules (in other words, if $\calI$ is a cogenerator),
then any left $\calI$-approximation is a monomorphism, and any $\calI$-coresolution is an exact sequence.
\end{rmk}

In the following, we only explain right $\calI$-approximations and $\calI$-resolutions, and omit the dual version.

By considering the case of $\calI = \add A$, it is easy to see that
  the concept of a $\calI$-resolution of $M$ is
  a generalization of the idea of a (not necessarily minimal) projective resolution of $M$.
  For minimality, we need the following definitions.

\begin{dfn}
  Let $M \in \mod A$.
  \begin{enumerate}
  \item
    The comma category $(\mod A)\comcat M$ from $\mod A$ to $M$ is defined as follows.
    Its objects are the morphisms $f \colon X \to M$ in $\mod A$,
    and for any objects $f\colon X \to M$, $f' \colon X' \to M$,
    the morphism set $\Hom(f, f')$ is given as the set of morphisms
    $g \colon X \to X'$ in $\mod A$ satisfying $f = f'g$.
    The composition is given by that of morphisms in $\mod A$.

    Note that for each $g \in \Hom(f, f')$,
    $g$ is an isomorphism in this category if and only if it is an isomorphism in $\mod A$.

  \item
    A morphism $f \colon X \to M$ in $\mod A$ is called \emph{right minimal}
    if every morphism $g \colon X \to X$ with $f = fg$
    (i.e., every $g \in \Hom(f,f)$) is an automorphism.

  \item
    Let $f, f'$ be morphisms in $(\mod A)\comcat M$.
    They are said to be \emph{equivalent} in $(\mod A)\comcat M$ if $\Hom(f,f') \ne \emptyset$
    and $\Hom(f',f) \ne \emptyset$.
    This is an equivalence relation on the set of morphisms in $(\mod A)\comcat M$.
    If $f'$ is right minimal and is equivalent to $f$, then $f'$ is called a {\em right minimal version} of $f$.
  \end{enumerate}
\end{dfn}

We give the following two properties of right minimal morphisms, which we will use later.
\begin{thm}[{\cite[Theorem 2.2]{auslander1997representation}}]
\label{thm:r-min}
  Let $f \colon X \to M$ in $\mod A$.
  Then there exists a decomposition $X = X' \ds X''$
  such that $f|_{X'}$ is a right minimal version of $f$ and $f|_{X''} = 0$.
  Moreover, $f|_{X'}$ is uniquely determined by $f$ up to equivalences.
\end{thm}

The following is immediate from the theorem above.

\begin{cor}
\label{cor:r-min-sect}
  Let $f \colon X \to M$ be a morphism in $\mod A$.
  Then the following are equivalent.
  \begin{enumerate}
  \item
    $f$ is right minimal.
  \item
    For any section $s \colon S \to X$, $fs=0$ implies $s=0$.
  \end{enumerate}
\end{cor}

The following can be shown using the corollary above, but here we give a simple alternative  proof.

\begin{lem}
\label{lem:ds-min}
  Let $f \colon X \to M$ and $f' \colon X' \to M'$ be morphisms in $\mod A$.
  If both of them are right minimal, then so is $f \ds f'$.
\end{lem}

\begin{proof}
  Let $T:=\spmat{a&b\\c&d} \in \End{(X\ds X')}$, $F:=f \ds f'=\spmat{f&0\\0&f'}$, and assume $FT=F$.
  This immediately yields $fa=f, f'c=0$, and $f'd=f'$,
  and we have $a \in \Aut{X}$ because $f$ is right minimal.
  Since $T=\spmat{\id&0\\ca^{-1}&\id}\spmat{a&b\\0&d-ca^{-1}b}$,
  it is sufficient to show that $u:=d-ca^{-1}b\in \Aut{X'}$.
  But it is clear from the right minimality of $f'$ and the fact that $f'u=f'd-f'ca^{-1}b=f'-0=f'$.
\end{proof}

Throughout the rest of this section, we let $G$ be a generator in $\mod A$,
and set $\calI:= \add G$, the full subcategory of $\mod A$ consisting of all direct summands of finite
direct sums of copies of $G$.
Then we have $A \in \calI$, and Remark~\ref{rmk:generator-ex} is applicable in this setting.

  Now we are ready to give the definition of minimal $\calI$-resolutions, which generalizes the definition of minimal projective resolutions.
\begin{dfn}
  Let $M$ be in $\mod A$.
  \begin{enumerate}
  \item
    A \emph{right minimal $\calI$-approximation} of $M$ is a right $\calI$-approximation
    $f \colon X \to M$ that is right minimal.
  \item
    A \emph{minimal $\calI$-resolution} of $M$ is
    a sequence
    \[
      \cdots \to X_2 \ya{f_2} X_1 \ya{f_1} X_0 \ya{f_0} M \to 0
    \]
    such that
    $f_0$ is a right minimal $\calI$-approximation of $M$, and for each
    $i \ge 1$, $f_i$ is a right minimal $\calI$-approximation of $\Ker f_{i-1}$.
    By Remark \ref{rmk:generator-ex}, this is automatically an exact sequence.
  \end{enumerate}
\end{dfn}

Theorem \ref{thm:r-min} gives the following link between minimal right $\calI$-approximations and not-necessarily-minimal ones.

\begin{lem}
\label{lem:min-nonmin-appx}
  Let $M$ be in $\mod A$, $f \colon I \to M$ a minimal right $\calI$-approximation, and
  $g \colon J \to M$ a right $\calI$-approximation.
  Then there exist an $I' \in \calI$ and an isomorphism $h \colon J \to I \ds I'$ such that
  $g = (f,0)h$, namely we have the commutative diagram
  \[
    \begin{tikzcd}
      J & M\\
      I \ds I' & M
      \Ar{1-1}{1-2}{"g"}
      \Ar{2-1}{2-2}{"{(f,0)}"'}
      \Ar{1-1}{2-1}{"h"}
      \Ar{1-2}{2-2}{equal}
  \end{tikzcd}.
\]
\end{lem}

Using Lemmas \ref{lem:ds-min} and \ref{lem:min-nonmin-appx}, we obtain the following that links
minimal resolutions and not-necessarily-minimal ones.

\begin{prp}
  Let $M$ be in $\mod A$,
  $I\down:= (\cdots \to I_1 \ya{f_1} I_0 \ya{f_0} M \to 0)$ a minimal $\calI$-resolution, and
  $J\down:= (\cdots \to J_1 \ya{g_1} J_0 \ya{g_0} M \to 0)$
  an $\calI$-resolution.
  Set $K_i:= \Im f_i$ and $L_i:= \Im g_i$ for all $i \ge 1$.

  Then for each $i \ge 0$, there exists $I'_i \in \calI$  such that
  $J\down \iso I\down \ds \cone(\id_{I'_0}) \ds \cone(\id_{I'_1})[1] \ds \cdots \ds \cone(\id_{I'_i})[i] \ds \cdots$
  and $L_i \iso K_i \ds I'_{i-1}$  for all $i \ge 1$.
  In particular, a minimal $\calI$-resolution is unique
  up to isomorphism.
\end{prp}

\begin{dfn}
  \label{dfn:dimension}
  Let $M$ be in $\mod A$.
  \begin{enumerate}
  \item If there exists an $\calI$-resolution $J\down$
    of $M$ of the form
    $J\down = (0 \to J_n \ya{g_n}  \cdots \to J_1 \ya{g_1} J_0 \ya{g_0} M \to 0)$
    for some $n \ge 0$ (noting that we placed no restrictions on $J_n,\hdots,J_1$ being nonzero),
    then we say that \emph{$\calI$-resolution dimension} of $M$ is at most $n$, and
    write $\calI\hyph\rdim M \le n$.
    Otherwise we say that $\calI$-resolution dimension of $M$ is infinity.

  \item If $\calI\hyph\rdim M \le n$ and $\calI\hyph\rdim M \not\le n-1$,
    then we say that $\calI$-resolution dimension of $M$ is equal to $n$, and denote it by
    $\calI\hyph\rdim M = n$.

  \item
    Finally, we set $\calI\hyph\rgldim A:= \sup\{\calI\hyph\rdim M \mid M \in \mod A\}$,
    and call it the \emph{$\calI$-resolution global dimension} of $A$.
  \end{enumerate}
\end{dfn}

\begin{lem}
  Let $M$ be in $\mod A$, $0 \le n$ an integer,
  $I\down = (I_i, f_i)_i$ a minimal $\calI$-resolution of $M$, and $J\down = (J_i, g_i)_i$ an $\calI$-resolution of $M$.
  Then the following are equivalent:
  \begin{enumerate}
  \item
    $\calI\hyph\rdim M \le n$.
  \item
    $\Im g_n \in \calI$.
  \item
    $I\down$ has the form $0 \to I_n \ya{f_n} \cdots \to I_1 \ya{f_1} I_0 \ya{f_0} M \to 0$.
  \end{enumerate}
\end{lem}

Set ${\La}:= \End_A(G)$, and regard $G$ as an $A$-${\La}$-bimodule\footnote{
Here $\End_A(G)$ denotes the algebra of ``right endomorphisms'' of $G$: $x(fg) = (xf)g$ for all $x \in G$ and $f,g \in \End_A(G)$.
Thus $\End_A(G) = (\Hom_A(G,G), \circ)\op$, where $\circ$ denotes the usual composition of maps: $(f\circ g)(x):= f(g(x))$ for all $x \in G$ and $f,g \in \Hom_A(G,G)$.}
Then we can consider the functor
$\Hom_A(G, \blank)\colon \mod A \to \mod {\La}$.
We denote by $\pd_{\La} Y$ the projective dimension of a left ${\La}$-module $Y$.
Now we give a way to compute the $\calI$-resolution dimension of an $A$-module
by the projective dimension of a ${\La}$-module.

The following is well-known (for instance, see \cite[Proposition~II.2.1]{auslander1997representation}
or \cite[Proposition~4.16]{blanchette2021homological}).
\begin{lem}
  \label{lem:Yoneda-like}
  Let $X, Y \in \mod A$ and $F \colon \Hom_A(G, X) \to \Hom_A(G, Y)$ be in $\mod \La$.
  If $X \in \add G$, then there exists a unique $f \in \Hom_A(X, Y)$ such that $F = \Hom_A(G, f)$.
\end{lem}

Using the above, the next Proposition follows.
(cf.\ \cite[Proposition~4.17(1)(2)]{blanchette2021homological})
\begin{prp}
\label{prp:I-rdim=pd}
Let $M$ be in $\mod A$.
Then
\[
\calI\hyph\rdim M = \pd_{\La} \Hom_A(G, M).
\]
\end{prp}

The following is immediate from Proposition \ref{prp:I-rdim=pd}
(cf.\ \cite[Proposition~4.17(3)]{blanchette2021homological}).
\begin{cor}
$\calI\hyph\rgldim A = \sup\{\pd_{\La} \Hom_A(G, M) \mid M \in \mod A\} \le \gldim {\La}$.
\end{cor}

In the case that $G$ is a generator-cogenerator, we get not only an upper bound for $\calI\hyph\rgldim A$ as above, but the following equality, which is immediate from Erdmann--Holm--Iyama--Schr{\"o}er \cite[Lemma 2.1]{erdmann2004radical}.
This will be useful for our setting.

\begin{prp}
\label{prp:intdim-gldim}
Assume further that $G$ is also a cogenerator,
i.e., all indecomposable injective $A$-modules are contained in $\calI$ up to isomorphisms.
Then 
\[
\calI\hyph\rgldim A = \gldim \La -2.
\]
\end{prp}

Finally, we recall the following definitions from \cite{blanchette2021homological}, in order to put our results in the context of their framework of homological invariants. For this part, we let $\calX$ be a finite set of indecomposable modules of $A$, $\calI = \add \calX$, and
let
\[
  \calE_\calX = \left\{0 \to L \ya{f_2} M \ya{f_1} N \to 0 \text{ short exact} \;\;\middle|\;\; f_1 \text{ is a right } \calI\text{-approximation} \right\}.
\]
\begin{dfn}[Relative Grothendieck groups]
  \leavevmode
  \begin{enumerate}
  \item Let $F$ be the free abelian group generated by the symbols $[M]$ for each isomorphism class of $M \in \mod A$. Let $H_\calX$ be the subgroup generated by $[M] - [L] - [N]$ for $(0 \to L \ya{f_2} M \ya{f_1} N \to 0) \in \calE_\calX$. The \emph{Grothendieck group of $A$ relative to $\calX$} is
    \[
      K_0(A, \calX) := F/H_\calX.
    \]
  \item Letting $\calE_{\text{min}}$ be the class of all split short exact sequences,
    the \emph{split Grothendieck group of $A$}, denoted
    $
      K_0^{\text{split}}(A) 
    $,
    is defined as above but replacing $\calE_\calX$ with $\calE_{\text{min}}$.
  \end{enumerate}
\end{dfn}

\begin{dfn}[Homological invariants \cite{blanchette2021homological}]
  Let $\calX$ be a finite set of indecomposable modules of $A$ and $\calI = \add \calX$.
  \begin{enumerate}

  \item An \emph{invariant} is a surjective group homomorphism $p: K_0^{\text{split}}(A) \rightarrow \bbZ^n$ for some $n\in\bbN$.
  \item Two invariants $p$ and $q$ are said to be equivalent if and only if $\ker p = \ker q$.

  \item An invariant $p$ is said to be a
    \emph{homological invariant relative to $\calX$} if all of the following hold:
    \begin{enumerate}
    \item All of the indecomposable projectives are in $\calI$.
    \item For each $M \in \mod A$, $\calI\hyph\rdim M < \infty$.
    \item $p$ is equivalent to the invariant given by the canonical quotient map
      \[
        K_0^{\text{split}}(A) \ya{\pi} K_0(A, \calX) \cong \bbZ^{|\calX|}
      \]
      where the last isomorphism follows from \cite[Proposition~4.9]{blanchette2021homological}.
    \end{enumerate}

  \end{enumerate}
\end{dfn}

\section{Finiteness and computation of the interval resolution global dimension}

\subsection{Finiteness of the interval resolution global dimension}

Throughout this subsection, we set $A := \k\calP$ to be the
incidence algebra of a poset $\calP$,
and $\calI = \add G$ with
$G:= \Ds_{I \in \bbI} V_I$, where $\bbI$ is the set of all intervals in $\calP$.
Where convenient, we abuse the notation and
  write $\bbI$ to also mean the set of modules $\{V_I \mid I \in \bbI\}$.
  With this abuse of notation, $\calI = \add \bbI$.
Modules in $\calI$ are called \emph{interval-decomposable} modules.

Since all projective indecomposable modules and all injective indecomposable modules are isomorphic to interval modules,
the module $G$ is indeed a generator and a cogenerator.

We specialize Definitions~\ref{dfn:approx}~and~\ref{dfn:dimension} to this setting.
\begin{dfn}
  Let $M$ be in $\mod \k\calP$ and $\calI = \add G$ with $G:= \Ds_{I \in \bbI(\calP)} V_I$.
  \begin{enumerate}
  \item
    A \emph{right interval-approximation} of $M$ is
    simply a right $\calI$-approximation of $M$.
  \item
    An {\em interval resolution} of $M$ is simply an $\calI$-resolution of $M$.
  \item The \emph{interval resolution dimension} of $M$ is simply the $\calI$-resolution dimension of $M$, and is denoted by
    $\intdim M$.
  \item
    The \emph{interval resolution global dimension} of $A$ is simply the $\calI$-resolution global dimension, and is denoted by $\intgldim(\k\calP)$.
  \end{enumerate}
\end{dfn}

In this subsection, we will prove that the global dimension of $\La = \End_A (G)$ is finite.
Then by Proposition \ref{prp:intdim-gldim}, we see that the interval resolution global dimension of $A$ is finite.
We use a Theorem in Ringel's work \cite{ringel2010iyama} for this.
Ringel~\cite{ringel2010iyama} gave a concise proof of Iyama’s finiteness theorem \cite{iyama2003finiteness} for the representation dimension by using the concept of a {\it left quasi-hereditary algebra}.
More precisely, he made use of a special type of that algebra with the additional property that all of its standard left modules have projective dimension at most $1$, which is called a
{\em left strongly quasi-hereditary algebra}.
Note that any left quasi-hereditary algebra has finite global dimension.

We cite the following from \cite[Theorem in \S 5]{ringel2010iyama} (cf. \cite[Lemma~2.2]{iyama2003finiteness}).
\begin{thm}
\label{thm:ringel-qh}
Let $B$ be an artin algebra, and $X$ a finitely generated $B$-module.
Then there exists a $B$-module $Y$ such that $C:= \End_B(X \oplus Y)$ is left strongly quasi-hereditary.
In particular, the global dimension of $C$ is finite.
Moreover, we can take this $Y$ as a module with the property that any indecomposable direct summand of $Y$ is a submodule of an indecomposable direct summand of $X$.
\end{thm}
For our purposes, we extract the following Corollary.
\begin{cor}
\label{cor:gen-statement}
Let $B$ be an artin algebra, and $X$ a finitely generated $B$-module.
Assume that for each indecomposable direct summand $X'$ of $X$,
all submodules of $X'$ are in $\add X$,
then $\End_B(X)$ is left strongly quasi-hereditary,
and its global dimension is finite.
\end{cor}

\begin{proof}
By Theorem \ref{thm:ringel-qh}, there exists a $B$-module $Y$ such that
$Y = \Ds_{i=1}^n Y_i$ with $Y_i \le X_i$ for some indecompsable
direct summand $X_i$ of $X$ for all $i = 1,\dots, n$
and that $\End_B(X \ds Y)$ is left strongly quasi-hereditary.
By assumption, all $Y_i$ are in $\add X$, and hence $Y \in \add X$.
Therefore, $\End_B(X \ds Y)$ is Morita equivalent to $\End_B(X)$.
Thus $\End_B(X)$ is also left strongly quasi-hereditary.
\end{proof}

\begin{lem}
\label{lem:intervalsub1}
Let $M$ be an interval module of $A=\k\calP$, and $N$ a submodule of $M$. Then, $N$ is interval-decomposable.
\end{lem}

\begin{proof}
We set $Q$ to be the Hasse diagram of $\calP$ regarded as a quiver.
Denote by $1_x \in \k$ the basis of $e_x N$
for all $x \in \calP$, where $e_x$ is the idempotent in $A$ associated to $x$.

We first show that the support
$\supp N:= \{x \in \calP \mid e_x N \ne 0\}$ of $N$ is convex.
Let $x, y \in \supp N$, $p$ a path from $x$ to $y$ in $Q$, and $z$ any vertex occurring in $p$.
It is enough to show that $e_z N \ne 0$.
Let $q$ be the subpath of $p$ from $x$ to $z$.
Since $M$ is an interval module, and $x,y \in \supp N \subseteq \supp M$, all the vertices  occurring in $p$ (and hence in $q$) are contained in $\supp M$.
Then since $N$ is a submodule of $M$, we have
\[
1_z = e_z q e_x 1_x \in e_z A e_x 1_x \le e_z N.
\]
Hence $e_z N \ne 0$, as desired.

Let $S_1, \dots, S_n$ be the connected components of $\supp N$.
Then for each $i = 1, \dots, n$, $S_i$ is  convex and connected, thus an interval.
Therefore, 
$N$ is the direct sum of the interval modules defined by the $S_i$'s.
\end{proof}

Corollary \ref{cor:gen-statement} and Lemma \ref{lem:intervalsub1} immediately imply
our desired result:

\begin{prp}
  \label{prp:intgldim-fin}
  Let $G:= \displaystyle\Ds_{I \in \bbI(\calP)} V_I$.
The global dimension of $\La := \End(G)$ is finite.
  Hence $\intgldim A$ is finite.
\end{prp}

\subsection{Computation of the interval resolution global dimension}

In this subsection, we give a way to compute $\intgldim A$ for an incidence algebra $A$
of a finite poset
by a computer. 
For this we first give a general way to compute $\calI$-resolution global dimension of
a finite-dimensional algebra in the following setting.

\begin{prp}
\label{prp:I-res-gldim}
Let $A$ be a finite-dimensional algebra,
$\bbJ$ a finite set of indecomposable $A$-modules,
$G:= \Ds_{M \in \bbJ}M$,
and $\calI:= \add G$.
Assume that $G$ is a generator and a cogenerator
and that $\calI\hyph\rgldim A < \infty$. Then
\[
\calI\hyph\rgldim A = \max_{M \in \bbJ} \calI\hyph\rdim \ta M
= \max_{M \in \bbJ} \calI\hyph\rdim \ta\inv M.
\]
\end{prp}

Our proof uses the following fundamental formula.

\begin{lem}
\label{lem:pd}
Let $\La$ be a finite-dimensional algebra.
If $0\to X \to Y \to Z \to 0$ is a short exact sequence of $\La$-modules, then we have
\[
\pd_\La Z \le \max\{\pd_\La Y, \pd_\La X +1\}.
\]
In particular, if $Y$ is projective, then
$\pd_\La Z \le \pd_\La X +1$.
\end{lem}

\begin{proof}
The short exact sequence above yields an exact sequence
\[
\Ext_\La^i(X,\hyph) \to \Ext_\La^{i+1}(Z,\hyph) \to \Ext_\La^{i+1}(Y,\hyph)
\]
of functors for all integers $i \ge 1$, which shows that
$\pd X \le i-1$ and $\pd Y \le i$ implies $\pd Z \le i$.
Hence the assertion follows by setting $i:= \max\{\pd_\La Y, \pd_\La X +1\} \ge 1$.
\end{proof}

\subsubsection*{Proof of Proposition \ref{prp:I-res-gldim}.}
Let $g \colon E_M \to M$ be a sink map.
Thus if $M$ is projective
$g \colon E_M = \rad M \to M$ is the inclusion and $\tau M = 0$,
or otherwise we have an almost split sequence
\[
0 \to\ta M \ya{f} E_M \ya{g} M \to 0.
\]
In either case, we have an exact sequence of left $\La$-modules of the form
\[
0 \to \Hom_A(G, \ta M) \to \Hom_A(G, E_M) \to \rad(G, M) \to 0.
\]
By Lemma \ref{lem:pd}, this shows
\begin{equation}
\label{eq:pd-rad}
\pd_\La \rad(G, M) \le \max\{\pd_\La \Hom_A(G, E_M), \pd_\La\Hom_A(G, \ta M) +1\}.
\end{equation}


Hence we have
\[
\begin{aligned}
\calI\hyph\rgldim A &=\gldim \La - 2 \quad \text{(by Proposition \ref{prp:intdim-gldim})}\\
&= \max\{\pd_\La(\Hom_A(G, M)/\rad(G, M)) \mid M \in \bbJ \}-2 \\
&\le \max\{\pd_\La \rad(G, M) + 1 \mid M \in \bbJ\}-2
\quad \text{(by Lemma \ref{lem:pd}, $Y:=\Hom_A(G, M)$ projective)}\\
 &= \max\{\pd_\La \rad(G, M) - 1 \mid M \in \bbJ\}\\
&\le \max\{\pd_\La \Hom_A(G, E_M) -1, \pd_\La\Hom_A(G, \ta M) \mid M \in \bbJ\}\ \text{(by \eqref{eq:pd-rad})}\\
&= \max\{\calI\hyph\rdim E_M -1, \calI\hyph\rdim \ta M \mid M \in \bbJ\} 
\ \text{(by Lemma \ref{prp:I-rdim=pd})}\\
&\le \calI\hyph\rgldim A.
\end{aligned}
\]
Therefore,
\[
\calI\hyph\rgldim A = \max\{\calI\hyph\rdim E_M -1, \calI\hyph\rdim \ta M \mid M \in \bbJ\}.
\]
By assumption $\calI\hyph\rgldim A = d$ for some non-negative integer $d$, and hence
there exists some $M \in \bbJ$ such that either $\calI\hyph\rdim E_M - 1 = d$ or $\calI\hyph\rdim \ta M = d$.
But the former is impossible because if this is the case, then
we have $\calI\hyph\rdim E_M = d + 1 > \calI\hyph\rgldim A = \max\{\calI\hyph\rdim X \mid X \in \mod A\}$, a contradiction.
As a consequence, 
$d = \max_{M \in \bbJ} \calI\hyph\rdim \ta M.$

The remaining equality is proved dually.
\qed
\medskip

\begin{rmk}
Since $\ta M = 0$ if $M$ is projective, and $\ta\inv M = 0$ if $M$ is injective,
we have
\[
\begin{aligned}
\max_{M \in \bbJ} \calI\hyph\rdim \ta M &=
\max\{0, \calI\hyph\rdim \ta M \mid M \in \bbJ, M \text{ is non-projective}\},\\
\max_{M \in \bbJ} \calI\hyph\rdim \ta\inv M &=
\max\{0, \calI\hyph\rdim \ta\inv M \mid M \in \bbJ, M \text{ is non-injective}\}.
\end{aligned}
\]
\end{rmk}

We apply Proposition~\ref{prp:I-res-gldim} to the setting of $A = \k\calP$ and
$\bbJ:= \{V_I \mid I \in \bbI\}$.
We note that
the assumption needed for Proposition~\ref{prp:I-res-gldim}
is guaranteed by Proposition~\ref{prp:intgldim-fin}.
Then we obtain the following.

\begin{prp}
\label{prp:intdim-ta}
Let $A = \k\calP$. Then we have
\[
\begin{aligned}
\intgldim A &= \max\{0, \intdim \ta V_I \mid I \in \bbI \text { with $V_I$ non-projective}\} = \max_{I \in \bbI}\intdim \ta V_I\\
&= \max\{0, \intdim \ta\inv V_I \mid I \in \bbI \text { with $V_I$ non-injective}\}
= \max_{I \in \bbI}\intdim \ta\inv V_I.
\end{aligned}
\]
\end{prp}


Using Proposition~\ref{prp:intdim-ta}
we computed several cases using the excellent GAP \cite{GAP4} package QPA \cite{qpa}.

\begin{exm}
Let $A = \k \Gf{m,n}$ ($m,n \ge 2$) and $\k = \mathbb{F}_2$, the finite field with $2$ elements.
Then the filled-in parts of the table below give the value (or a lower bound)
obtained by numerical computation of $\intgldim \k \Gf{m,n}$
in the row labelled $m$ and column labelled $n$.
By definition of $A = \k \Gf{m,n}$,
the table is symmetric so that it suffices to compute only the upper or lower part.
\begin{center}
\begin{tabular}{ l|ccccccccc } 
   & 2 & 3 & 4 & 5 & 6 & 7 & 8 & 9 & 10 \\ 
  \hline
 2 & 0 & 1 & 2 & 2 & 2 & 2 & 2 & 2 & 2\\ 
 3 & 1 & 2 & 3 & 4 & 4 & 4 \\
 4 & 2 & 3 & 4 & 5 & $\geq$\,6 & $\geq$\,6 \\
 5 & 2 & 4 & 5 \\
 6 & 2 & 4 & $\geq$\,6 \\
 7 & 2 & 4 & $\geq$\,6 \\
\end{tabular}
\end{center}
\end{exm}

These computations suggest some conjectures on the values
of $\intgldim \k \Gf{m,n}$ as follows. The first conjecture is about the commutative ladders \cite{escolar2016persistence}.
\begin{cnj}
\label{cnj:ladder}
For the algebra $A = \k \Gf{m,n}$,
$\intgldim A = 2$ if $n=2$ and $m \geq 4$.
\end{cnj}

If the next conjecture is true, then the preceding conjecture also holds with $C(2) = 2$.
\begin{cnj}
\label{cnj:stabilize}
There exists a function $C : \bbN \rightarrow \bbN \cup \{0\}$ such that 
for each $m \in \bbN$,
\[
\intgldim \k \Gf{m,n} = C(m)
\]
for all $n \geq m+2$.
That is, considering $\intgldim \k \Gf{m,n}$ as a function of $n$ for fixed $m$, this function stabilizes to some constant $C(m)$ starting from $m+2$.
\end{cnj}

\section{Relationship between interval-decomposable approximations and interval resolutions for commutative ladders}

In this section, we relate the interval resolutions of persistence modules in the
$2\times n$ case with a modified version of the compressed multiplicity invariant \cite{asashiba2019approximation}.
Via M\"obius inversion, we obtain an invariant that we show to be a homological invariant in Corollary~\ref{cor:homologicalinvariant}.

To fix some notation, we let $\Gf{m,2} = (Q, \ro)$ be a commutative ladder, where
$Q$ is presented as follows:
\[
  \begin{tikzcd}
    y_1 & y_2 & \cdots & y_n\\
    x_1 & x_2 & \cdots & x_n
    \Ar{1-1}{1-2}{"c_1"}
    \Ar{1-2}{1-3}{"c_2"}
    \Ar{1-3}{1-4}{"c_{n-1}"}
    \Ar{2-1}{2-2}{"a_1"'}
    \Ar{2-2}{2-3}{"a_2"'}
    \Ar{2-3}{2-4}{"a_{n-1}"'}
    \Ar{2-1}{1-1}{"b_1"}
    \Ar{2-2}{1-2}{"b_2"'}
    \Ar{2-3}{1-3}{"\cdots", phantom}
    \Ar{2-4}{1-4}{"b_n"'}
  \end{tikzcd}.
\]
and $\ro$ the full commutativity relations.
Note that compared to Diagram~\eqref{eq:commgrid}, we have
set $x_i:= (i,1)$ and $y_i:= (i,2)$,
$b_i:= (1,\be_i)$ for all $i = 1,\dots, n$, and
$a_i:= (\al_i, 1), c_i:= (\al_i, 2)$ for all for all $i = 1,\dots, n-1$
to make the notation simpler.
We regard the algebra $A:= \k(Q,\ro)$ as a category defined by this quiver with
the full commutativity relations $\ro$ as in the general way.
Namely, in general, the path algebra $\k(\Ga,R)$ of a bound quiver $(\Ga, R)$ is regarded as the category $\calC$, where 
the set of objects of $\calC$ is given by $\Ga_0$, for any $x, y \in \Ga_0$,
$\calC(x,y):= e_y\k(\Ga, R)e_x$, and the composition of $\calC$ is given by the
multiplication of $\k(\Ga, R)$.
This $\calC$ is called the {\em path category} of the bound quiver $(\Ga,R)$.

\subsection{Compression}
First, we introduce a modified version of compression of modules that was studied in \cite{asashiba2019approximation}.

Let us fix some notation.
Let $I = [x_i, x_j]_1 \sqcup [y_k, y_l]_2$ be an interval of $Q$
(recall the the discussion in Section~4.1 of \cite{asashiba2022interval}
concerning the ``staircase'' shape of intervals).
We have $1 \le k \le i \le l \le j \le n$, and $I$ is
illustrated as follows:
\begin{equation}
\label{eq:fixed-interval}
\begin{tikzcd}
y_k & \cdots & y_i & \cdots & y_l\\
    &        & x_i & \cdots & x_l &\cdots & x_j
\Ar{1-1}{1-2}{"c_k"}
\Ar{1-2}{1-3}{"c_{i-1}"}
\Ar{1-3}{1-4}{"c_{i}"}
\Ar{1-4}{1-5}{"c_{l-1}"}
\Ar{2-3}{2-4}{"a_i"'}
\Ar{2-4}{2-5}{"a_{l-1}"'}
\Ar{2-5}{2-6}{"a_{l}"'}
\Ar{2-6}{2-7}{"a_{j-1}"'}
\Ar{2-3}{1-3}{"b_i"}
\Ar{2-4}{1-4}{"\cdots", phantom}
\Ar{2-5}{1-5}{"b_l"'}
\end{tikzcd}
\end{equation}

On the other hand, let $Q'$ be the quiver:
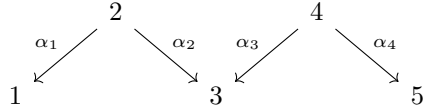
\[
  \begin{tikzcd}
    & 2 && 4\\
    1 && 3 && 5
    \Ar{1-2}{2-1}{"\al_1"'}
    \Ar{1-2}{2-3}{"\al_2"}
    \Ar{1-4}{2-3}{"\al_3"'}
    \Ar{1-4}{2-5}{"\al_4"}
  \end{tikzcd}
\]
and set $B:= \k Q'$ to be the path category of the quiver $Q'$.

Moreover, for $i,j$ in $Q'_0$, we set $\intval{i,j}:= \{x \in Q'_0 \mid i \le x \le j\}$, where $\le$ is the total order in $Q'_0:=\{1,2,3,4,5\}$ as a subset of the integers.
Note that $\intval{i,j}$ is an interval of the poset\footnote{
Here, $Q'$ is regarded as the poset $(Q'_0, \preceq)$ with the partial order $\preceq$
defined by $x \preceq y$ iff there exists an arrow $x \to y$ for all $x,y \in Q'_0$.
} $Q'$.
For each vertex $x$ of $Q$ (resp.\ $Q'$) we denote by $e_x$ (resp.\ $e'_x$)
the path of length 0 at $x$, and for each morphism $f$ in $\k Q$,
$\ovl{f}$ denotes the image of $f$ under the canonical morphism $\k Q \to A$.

For the interval $I$, we define a quiver morphism $\xi_I \colon Q' \to A$ (depending on the form of $I$).
\begin{itemize}
\item If $I$ is of the form \eqref{eq:fixed-interval}, that is, if $I$ spans two rows of the commutative ladder,
  define $\xi_I \colon Q' \to A$ by the following table:
  \begin{center}
    \begin{tabular}{c||c|c|c|c|c|c|c|c|c}
      $x$ &  1 & 2 & 3 & 4 & 5 & $\al_1$ & $\al_2$ & $\al_3$ & $\al_4$\\
      \hline
      $\xi_I(x)$ & $y_l$ & $y_k$ & $y_i$ & $x_i$ & $x_j$& $\ovl{c_{l-1}\cdots c_k}$ &$\ovl{c_{i-1}\cdots c_k}$ & $\ovl{b_i}$ & $\ovl{a_{j-1}\cdots a_i}$
    \end{tabular},
  \end{center}
  where if $k=i$, then $c_{i-1}\cdots c_k$ is replaced by $e_i$,
  (similar for the cases $i=l$ or $l=j$).

\item When $I$ is of the form
  $x_i \to \cdots \to x_j$ with $0 \le i \le j \le n$
  ($I$ is contained in the lower row of the commutative ladder),
  we define $\xi_I$ by the following table:
  \begin{center}
    \begin{tabular}{c||c|c|c|c|c|c|c|c|c}
      $x$ &  1 & 2 & 3 & 4 & 5 & $\al_1$ & $\al_2$ & $\al_3$ & $\al_4$\\
      \hline
      $\xi_I(x)$ & $x_j$ & $x_i$ & $x_i$ & $x_i$ & $x_j$& $\ovl{a_{j-1}\cdots a_i}$ &$\ovl{e_i}$ & $\ovl{e_i}$ & $\ovl{a_{j-1}\cdots a_i}$
    \end{tabular}.
  \end{center}

\item We make a similar construction of $\xi_I$ when $I$ is of the form
  $y_k \to \cdots \to y_l$ with $0 \le k \le l \le n$
  ($I$ is contained in the upper row of the commutative ladder)
by the following table:
  \begin{center}
    \begin{tabular}{c||c|c|c|c|c|c|c|c|c}
      $x$ &  1 & 2 & 3 & 4 & 5 & $\al_1$ & $\al_2$ & $\al_3$ & $\al_4$\\
      \hline
      $\xi_I(x)$ & $y_l$ & $y_k$ & $y_k$ & $y_k$ & $y_l$& $\ovl{c_{l-1}\cdots c_k}$ &$\ovl{e_k}$ & $\ovl{e_k}$ & $\ovl{c_{l-1}\cdots c_k}$
    \end{tabular}.
  \end{center}
\end{itemize}

In the first case, $\xi_I$ is visualized as follows:
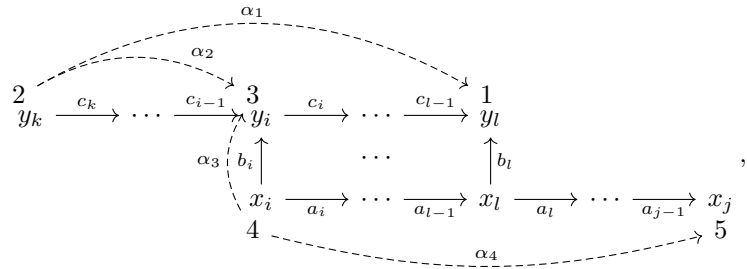
\[
\begin{tikzcd}
y_k & \cdots & y_i & \cdots & y_l\\
    &        & x_i & \cdots & x_l &\cdots & x_j
\Ar{1-1}{1-2}{"c_k"}
\Ar{1-2}{1-3}{"c_{i-1}"}
\Ar{1-3}{1-4}{"c_{i}"}
\Ar{1-4}{1-5}{"c_{l-1}"}
\Ar{2-3}{2-4}{"a_i"'}
\Ar{2-4}{2-5}{"a_{l-1}"'}
\Ar{2-5}{2-6}{"a_{l}"'}
\Ar{2-6}{2-7}{"a_{j-1}"'}
\Ar{2-3}{1-3}{"b_i"}
\Ar{2-4}{1-4}{"\cdots", phantom}
\Ar{2-5}{1-5}{"b_l"'}
\end{tikzcd}
\hspace{-98mm}
\begin{tikzcd}[column sep=13mm, row sep=13mm]
2 &  & 3 &  & 1\\
    &        & 4 &  &  & & 5
\Ar{1-1}{1-5}{bend left, dashed,"\alpha_1"}
\Ar{1-1}{1-3}{bend left, dashed,"\alpha_2", near end}
\Ar{2-3}{1-3}{bend left, dashed,"\alpha_3"}
\Ar{2-3}{2-7}{bend right=5mm, dashed,"\alpha_4"}
\end{tikzcd},
\]
where each broken arrow $\alpha$ represents an arrow in the quiver $Q'$,
and the corresponding solid path represents its image $\xi_I(\alpha)$ in the path category $A$.

In any case, $\xi_I$ uniquely extends to a linear functor $F_I \colon B \to A$.
By using $F_I$, we regard $A$ to be the
$A$-$B$-bimodule $\bi{A}AB = A(F_I(\cdot),\blank)$\footnote{Note that $A$ is a path category. Here, we refrain from writing $\Hom_A(F_I(\cdot),\blank)$ because usually $\Hom_A(X,Y)$ is used for $X, Y$ $A$-modules, but here in $A(F_I(x),y)$, we mean the set of morphisms from $F_I(x)$ to $y$ in the category
  $A$.},
which gives us an adjoint pair
\[
  \begin{tikzcd}
    \mod A & \mod B
    \Ar{1-1}{1-2}{"R_I"', bend right=20}
    \Ar{1-2}{1-1}{"L_I"', bend right=20}
  \end{tikzcd},
\]
where $\mod A$ (resp.\ $\mod B$) is the category of the finite-dimensional (left) $A$-modules (resp.\ $B$-modules), and
$L_{I}:= \bi{A}AB \otimes_B\blank$ is a left adjoint to
$R_{I}:= \Hom_A(\bi{A}AB, \blank)$.
Recall that an $A$-module $M$ is a functor $A \to \mod \k$, and then
$R_I(M) \iso M \circ F_I$.
For instance, as is easily seen, $R_I(V_I) \iso V_{\intval{1,5}}$
for all $I \in \bbI$.

\begin{dfn}[Compressed multiplicity]
  \label{def:compressedmult}

  Let $\xi$ be the function associating $I \in \bbI$ to $\xi_I$ as defined above.
  We define the \emph{compressed multiplicity with respect to $\xi$} of $V_I$ in $M$ as
  \[
    c_M^{\xi}(I):= d_{R_{I}(M)}(R_{I}(V_I)),
  \]
  which is multiplicity of $R_{I}(V_I)$ as a direct summand of $R_{I}(M)$,
  for all $I \in \bbI$.
\end{dfn}
This is a modification of the compressed multiplicity introduced in \cite{asashiba2019approximation}.
In that previous work, for $F_I$ we used a functor defined
by the inclusion of essential vertices, and so the corresponding
$R_I(M) \iso M \circ F_I$ is simply the compression functor of \cite{asashiba2019approximation}.
As also noted in \cite{asashiba2019approximation}, when $F_I$ is defined using the inclusion of $I$ as is, the corresponding compressed multiplicity is equal to the generalized rank invariant of Kim and Memoli~\cite{kim2018generalized} (see \cite{asashiba2019approximation} for a more detailed discussion).
  Here, we emphasize that instead of inclusion,
  $F_I$ is obtained from the quiver morphism $\xi_I$ as defined above, thus yielding a new variant of the invariant.

\begin{exm}
Let us compare the compressed multiplicity $c_M^{\xi}$ with the compressed multiplicities $c_M^{ss}$ and $c_M^{cc}$ introduced in \cite{asashiba2019approximation} (refer to \cite{asashiba2019approximation} for the definitions).
We set
\[
M:=\begin{tikzcd}
\k & \k^2 & \k & 0\\
0  & \k & \k^2 & \k 
\Ar{1-1}{1-2}{"\spmat{1\\1}"}
\Ar{1-2}{1-3}{"(0\;\;1)"}
\Ar{1-3}{1-4}{"0"}
\Ar{2-1}{2-2}{"0"'}
\Ar{2-2}{2-3}{"\spmat{1\\0}"'}
\Ar{2-3}{2-4}{"(1\;\;1)"'}
\Ar{2-1}{1-1}{"0"'}
\Ar{2-2}{1-2}{"\spmat{0\\1}"}
\Ar{2-3}{1-3}{"(1\;\;0)"}
\Ar{2-4}{1-4}{"0"}
\end{tikzcd} 
\text{as a representation of $\Gf{4,2}=$}
\begin{tikzcd}
y_1 & y_2 & y_3 & y_4\\
x_1 & x_2 & x_3 & x_4
\Ar{1-1}{1-2}{"c_1"}
\Ar{1-2}{1-3}{"c_2"}
\Ar{1-3}{1-4}{"c_3"}
\Ar{2-1}{2-2}{"a_1"'}
\Ar{2-2}{2-3}{"a_2"'}
\Ar{2-3}{2-4}{"a_3"'}
\Ar{2-1}{1-1}{"b_1"'}
\Ar{2-2}{1-2}{"b_2"}
\Ar{2-3}{1-3}{"b_3"}
\Ar{2-4}{1-4}{"b_4"}
\end{tikzcd},
\]
and let $I := [x_2, x_3]_1 \sqcup [y_1, y_3]_2$
and $J := [x_2, x_4]_1 \sqcup [y_2, y_3]_2$. Then we have
$c_M^{\xi}(I)=0$ and $c_M^{\xi}(J)=1$.
On the other hand, it is easy to verify that $c_M^{ss}(I)=1$ and $c_M^{cc}(J)=0$, and thus
$c_M^{ss} \neq c_M^{\xi} \neq c_M^{cc}$.
Furthermore $c_M^{cc}(I)=0$, and thus
we also see that $c_M^{ss} \neq c_M^{cc}$.
\end{exm}

Let us relate the above compressed multiplicities to the interval resolutions studied in the previous sections.
First, we lay down some groundwork.

\begin{lem}
  For each $I\in\mathbb{I}$ and each $x \in Q'_0$, we have $L_I(Be'_x) \iso Ae_{\xi_I(x)}$.
\end{lem}
\begin{proof}
  Since the canonical morphism $\bi{A}AB \otimes_B \bi{B}BB \to \bi{A}AB$
  of $A$-$B$-bimodules defined by sending
  $a \otimes b$ to $a{F}_I(b)$ for all morphisms $a$ in $A$
  (resp.\ $b$ in $B$) is an isomorphism,
  we have $A\otimes_B Be'_x \iso AF_I(e'_x) = Ae_{\xi_I(x)}$.
\end{proof}

The functor $L_I: \mod B \rightarrow \mod A$ satisfies the following nice property.
\begin{prp}\label{prp:L-int-dec}
  For any $I\in\mathbb{I}$, and for any indecomposable module X in $\mod B$, $L_I(X)$ is isomorphic to an interval module in $\mod A$, or is zero.
  Therefore, $L_I(Y)$ is interval decomposable or zero for all modules $Y$ in $\mod B$.
\end{prp}
\begin{proof}
  $B$ has $15$ isoclasses of indecomposable modules, whose complete set of representatives
  is given by the interval modules $V_{\intval{i,j}}$ for all $i,j$ with $1 \le i \le j \le n$.
  The assertion is verified by explicit computation of $L_I(X)$ using
  a minimal projective presentation of $X$ and the lemma above.
  For instance, we check the statement for $X = V_{\intval{2,4}}$.
  This module has the following minimal projective presentation:
  \[
    0 \to Be'_1 \ds Be'_3 \ds Be'_5 \ya{\spmat{\al_1&\al_2&0\\0&\al_3&\al_4}} Be'_2 \ds Be'_4 \to X \to 0,
  \]
  where $\al_i$ stands for the right multiplication by $\al_i$ for each $i=1,\dots,4$.
  Then since $L_I$ is right exact, we have an exact sequence
  \[
    Ae_{y_l} \ds Ae_{y_i} \ds Ae_{x_j} \ya{\spmat{\al_1'&\al_2'&0\\0&\al_3'&\al_4'}} Ae_{y_k} \ds Ae_{x_i} \to L_I(X) \to 0,
  \]
  where $\al_i'$ stands for the morphism $A \otimes_B \al_i$ for each $i=1,\dots,4$.
  By computing the cokernel of the morphism, we have
  \[
    L_I(X) \iso V_{[x_i,x_{j-1}]_1 \sqcup [y_k,y_{l-1}]_2}.
  \]
  (We set $V_{\emptyset}:= 0$. Thus if $i=j$ and $k=l$, then $L_I(X)=0$.)

  Similarly, the other cases can be verified.
\end{proof}

We then recall that for $M$  in $\mod A$, we have the following.
\begin{enumerate}
\item
  Since $A = \k \Gf{m,2}$ has a finite interval resolution global dimension (by Proposition~\ref{prp:intgldim-fin}), say $r$, 
  there exists an interval resolution of $M$ in the form of the following exact sequence
  with some non-negative integers $d_J^{(i)}$ for all $i=0,1,\dots, r$ and $J \in \bbI$:
  \begin{equation}\label{eq:int-resol}
    0 \to \Ds_{J \in \bbI}V_J^{d_J^{(r)}}
    \ya{f_r} \cdots
    \ya{f_2}
    \Ds_{J \in \bbI}V_J^{d_J^{(1)}} \ya{f_1}
    \Ds_{J \in \bbI}V_J^{d_J^{(0)}} \ya{f_0} M \to 0.
  \end{equation}
\item
  If $X \in \calI$, then the sequence \eqref{eq:resol-ex}
  and the exactness of the sequence \eqref{eq:int-resol}
  yields an exact sequence
    \begin{equation}\label{eq:expl-resol-ex}
      \begin{aligned}
        0 \to [X,\Ds_{J \in \bbI}V_J^{d_J^{(r)}}] & \ya{[X,f_r]} \cdots
                                                    \ya{[X,f_2]}
        [X,\Ds_{J \in \bbI}V_J^{d_J^{(1)}}]\\
        &\ya{[X,f_1]}
          [X,\Ds_{J \in \bbI}V_J^{d_J^{(0)}}] \ya{[X,f_0]}
          [X,M] \to 0.
      \end{aligned}
    \end{equation}
    where $[\cdot,\blank]:=\Hom_A(\cdot,\blank)$ for short.
    Note that the requirement that $X \in \calI$ is needed here.
\end{enumerate}

With the preparations finished,
we give the following Theorem~\ref{thm:comp-mult-resol} relating the interval resolution \eqref{eq:int-resol}
with the compressed multiplicities $c^{\xi}_M$ with respect to $\xi$ of $M$.
\begin{thm}
\label{thm:comp-mult-resol}
  Let $M$ be in $\mod A$ with an interval resolution \eqref{eq:int-resol}.
  Then for any $I \in \mathbb{I}$, we have
  \[
    c^{\xi}_M(I) = \sum_{I \le J \in \bbI}\sum_{i=0}^r (-1)^i d_J^{(i)},
  \]
  where $I \le J$ means that $I$ is a subquiver of $J$.
\end{thm}

\begin{proof}
  Recall that $R_I(V_I) = V_{\intval{1,5}}$.
  Then the almost split sequence starting from $R_I(V_I)$ has the following form:
  \[
    0 \to R_I(V_I) \to V_{\intval{2,5}} \ds V_{\intval{1,4}} \to V_{\intval{2,4}} \to 0,
  \]
  which yields a minimal projective resolution
  \[
    0 \to \Hom_B(V_{\intval{2,4}},\blank) \to \Hom_B(V_{\intval{2,5}} \ds V_{\intval{1,4}},\blank)
    \to \Hom_B(R_I(V_I),\blank) \to \calS_{R_I(V_I)} \to 0,
  \]
  of the simple functor $\calS_{R_I(V_I)}:= \Hom_B(R_I(V_I),\blank)/\rad_B(R_I(V_I, \blank)$
  corresponding to $R_I(V_I)$ by Auslander-Reiten theory.
  This shows that
  \begin{equation}\label{eq:AR-RM}
    \begin{aligned}
      d_{R_I(M)}(R(V_I)) &= \dim\calS_{R_I(V_I)}(R_I(M))\\
      &=\dim\Hom_B(R_I(V_I),R_I(M)) -\dim\Hom_B(V_{\intval{2,5}} \ds V_{\intval{1,4}},R_I(M))\\
      &\phantom{=} {} + \dim\Hom_B(V_{\intval{2,4}},R_I(M))
    \end{aligned}
  \end{equation}


  Now since $\bi{A}AB$ is projective as a left $A$-module,
  $R_I = \Hom_A(\bi{A}AB, \blank)$ is an exact functor.
  Then the exact sequence~\eqref{eq:int-resol} gives us an exact sequence
  \begin{equation}\label{eq:R(int-resol)}
    0 \to \Ds_{J \in \bbI}R_I(V_J)^{d_J^{(r)}} \ya{R_I(f_r)}\cdots\ya{R_I(f_r)}
    \Ds_{J \in \bbI}R_I(V_J)^{d_J^{(1)}} \ya{R_I(f_1)}
    \Ds_{J \in \bbI}R_I(V_J)^{d_J^{(0)}} \ya{R_I(f_0)} R_I(M) \to 0.
  \end{equation}

  Let $Y$ be in $\mod B$. Apply $\Hom_B(Y,\blank)$ to the exact sequence~\eqref{eq:R(int-resol)},
  and consider the following commutative diagram obtained by the adjoint pair $(L_I, R_I)$:
  \[\small
    \begin{tikzcd}[column sep=15pt]
      0 & (Y,\DS_{J \in \bbI}R_I(V_J)^{d_J^{(r)}}) &\cdots&
      (Y,\DS_{J \in \bbI}R_I(V_J)^{d_J^{(1)}}) &
      (Y,\DS_{J \in \bbI}R_I(V_J)^{d_J^{(0)}}) &
      (Y,R_I(M)) & 0
      \\
      0 & {[L_I(Y),\DS_{J \in \bbI}V_J^{d_J^{(r)}}]} &\cdots&
      {[L_I(Y),\DS_{J \in \bbI}V_J^{d_J^{(1)}}]} &
      {[L_I(Y),\DS_{J \in \bbI}V_J^{d_J^{(0)}}]} &
      {[L_I(Y),M]} & 0,
      \Ar{1-1}{1-2}{}
      \Ar{1-2}{1-3}{}
      \Ar{1-3}{1-4}{}
      \Ar{1-4}{1-5}{}
      \Ar{1-5}{1-6}{}
      \Ar{1-6}{1-7}{}
      \Ar{2-1}{2-2}{}
      \Ar{2-2}{2-3}{}
      \Ar{2-3}{2-4}{}
      \Ar{2-4}{2-5}{}
      \Ar{2-5}{2-6}{}
      \Ar{2-6}{2-7}{}
      \Ar{1-2}{2-2}{}
      \Ar{1-4}{2-4}{}
      \Ar{1-5}{2-5}{}
      \Ar{1-6}{2-6}{}
    \end{tikzcd}
  \]
  where we set $(\cdot,\blank):=\Hom_B(\cdot,\blank)$ and
  $[\cdot,\blank]:=\Hom_A(\cdot,\blank)$ for short.
  By Proposition \ref{prp:L-int-dec}, $L_I(Y) \in \calI$ and thus
  the lower row is exact (replace $X$ by $L_I(Y)$ in the exact sequence \eqref{eq:expl-resol-ex}).
  Since all the vertical maps are isomorphisms, the upper row is also exact.

  This yields the following equality:
  \[
    \dim\Hom_B(Y,R_I(M)) = \sum_{J\in \bbI}\left(\sum_{i=0}^r (-1)^id_J^{(i)}\right)
    \dim\Hom_B(Y,R_I(V_J)).
  \]
  By applying this formula to Equation~\eqref{eq:AR-RM} for
  $Y = R_I(V_I)$, $Y = V_{\intval{2,5}} \ds V_{\intval{1,4}}$, 
  and $Y = V_{\intval{2,4}}$, we have
  \begin{equation}\label{eq:almost-done}
    \begin{aligned}
      d_{R_I(M)}(R_I(V_I))
      &= \sum_{J\in \bbI}\left(\sum_{i=0}^r (-1)^id_J^{(i)}\right)
      (\dim\Hom_B(R_I(V_I),R_I(V_J))\\
      &\hspace{1.5em}-\dim\Hom_B(V_{\intval{2,5}} \ds V_{\intval{1,4}},R_I(V_J))+\dim\Hom_B(V_{\intval{2,4}},R_I(V_J)))\\
      &= \sum_{J\in \bbI}\left(\sum_{i=0}^r (-1)^id_J^{(i)}\right)d_{R_I(V_J)}(R_I(V_I)),
    \end{aligned}
  \end{equation}
  where the last equality follows from Equation~\eqref{eq:AR-RM} with
  $M = V_J$.
  Now since every entry of the dimension vector of $R_I(V_J)$ is at most $1$,
  we have
  \[
    d_{R_I(V_J)}(R_I(V_I)) = 0 \text{ or } 1.
  \]
  If $I$ is a subquiver of $J$, then it is obvious that $R_I(V_J) = R_I(V_I)$, and
  hence $d_{R_I(V_J)}(R_I(V_I)) = 1$.
  Conversely, if $d_{R_I(V_J)}(R_I(V_I)) = 1$,
  then $\dim R_I(V_J)(x) = \dim R_I(V_I)(x) \ne 0$
  for all $x \in \{x_i, x_j, y_k, y_l\}$,
  and hence $\{x_i, x_j, y_k, y_l\} \subseteq J_0$.
  Therefore, $I_0 \subseteq J_0$, and $I$ is a subquiver of $J$.
  As a consequence,
  \[
    d_{R_I(V_J)}(R_I(V_I)) =
    \begin{cases}
      1 & \text{if } I \le J,\\
      0 & \text{otherwise}.
    \end{cases}
  \]
  This together with Equation~\eqref{eq:almost-done} proves the assertion.
\end{proof}

That is, the compressed multiplicity $c^{\xi}_M$ with respect to $\xi$ of $M$ can be expressed
in terms of a formula involving only the multiplicities of the intervals in an interval resolution of $M$.
Following the ideas in previous works \cite{kim2018generalized, asashiba2019approximation},
we use M\"obius inversion to obtain another invariant.
First, we note that $\mathbb{I}$ can be given the structure of a poset by setting $I \leq J$ if and only if $I \subseteq J$, for all $I, J \in \mathbb{I}$.
\begin{dfn}[Interval approximation]
Recall that
$
  c_M^{\xi}(I):= d_{R_{I}(M)}(R_{I}(V_I))
$
for all $I \in \bbI$.
We define \emph{interval approximation $\de_M^{\xi}$ with respect to $\xi$} to be the
the M\"obius inversion of $c_M^{\xi}$, which is defined by
\[
  \de_M^{\xi}(J):= \sum_{S\subseteq \cov(J)} (-1)^{\# S}c_M^{\xi}(\bigvee S)
\]
for all $J \in \bbI$.
\end{dfn}

By the general theory of M\"obius inversion, and since by Theorem~\ref{thm:comp-mult-resol}
\[
  c^{\xi}_M(I) = \sum_{I \le J \in \bbI} \left(\sum_{i=0}^r (-1)^i d_J^{(i)}\right),
\]
we obtain the following.
\begin{cor}
\label{cor:de-xi}
  Let $M$ be in $\mod A$ with an interval resolution \eqref{eq:int-resol}.
  Then we have
  \[
    \de_M^{\xi}(J) = \sum_{i=0}^r (-1)^id_J^{(i)}
  \]
  for all $J \in \bbI$.
\end{cor}

\begin{rmk}
  Assume that $M$ is interval-decomposable and that the interval resolution \eqref{eq:int-resol} is a minimal interval resolution.
  Then we have $d_J^{(i)} = 0$
  for all $i\ge 1$ and all $J \in \bbI$, and
  $M \iso \Ds_{J \in \bbI}V_J^{d_J^{(0)}}$.
  Therefore, $\de_M^{\xi}(J) = d_J^{(0)} = d_M(V_J)$ for all $J \in \bbI$.

  That is, for $M$ interval-decomposable, the value $\de_M^{\xi}(J)$ of interval approximation at $J$
  is exactly equal to the multiplicity of the interval module $V_J$ as a direct summand of $M$.
\end{rmk}

Furthermore, expressed in the framework of \cite{blanchette2021homological},
our results Proposition~\ref{prp:intgldim-fin} and Corollary~\ref{cor:de-xi} translates to the following Corollary~\ref{cor:homologicalinvariant}.

\begin{cor}
  \label{cor:homologicalinvariant}
  Let $A = \k \Gf{m,2}$, and let $M$ be in $\mod A$.
  Then, interval approximation with respect to $\xi$,
  which associates $M$ to $\de_M^{\xi} : \bbI \rightarrow \mathbb{Z}$,
  is the homological invariant relative to $\bbI$ (and thus equivalent to the dim-hom invariant relative to $\bbI$ by \cite[Theorem~1.1]{blanchette2021homological}).
\end{cor}
\begin{proof}
  Note first that
  $\bbI$ contains all indecomposable projectives, and that
  $A = \k \Gf{m,2}$ has a finite interval resolution global dimension by Proposition~\ref{prp:intgldim-fin}.

  For
  $M \in \mod A$ with the interval resolution
  \[
    0 \to \Ds_{J \in \bbI}V_J^{d_J^{(r)}}
    \ya{f_r} \cdots
    \ya{f_2}
    \Ds_{J \in \bbI}V_J^{d_J^{(1)}} \ya{f_1}
    \Ds_{J \in \bbI}V_J^{d_J^{(0)}} \ya{f_0} M \to 0,
  \]
  the equivalence class of $M$ in $K_0(A, \bbI)$ is 
  \[
    [M] = \sum_{i=0}^r (-1)^i \left[\Ds_{J \in \bbI}V_J^{d_J^{(i)}}\right]
  \]
  by the usual argument of breaking an exact sequence into short exact sequences.
  Thus,
  \[
    [M] = \sum_{J \in \bbI} \sum_{i=0}^r (-1)^i d_J^{(i)} [V_J] = \sum_{J \in \bbI} \de_M^{\xi}(J) [V_J]
  \]
  by applying Corollary~\ref{cor:de-xi}.
  Under the identifications $K_0(A, \bbI) \cong \bbZ^{|\bbI|} \cong \bbZ^{\bbI}$ (functions from $\bbI$ to $\bbZ$) using \cite[Proposition~4.9]{blanchette2021homological},  this simply means that  $[M] = \de_M^{\xi}$. Thus, the canonical quotient map
  \[
    \pi: K_0^{\text{split}}(A) \rightarrow K_0(A, \bbI)
  \]
  takes the equivalence class of $M \in \mod{A}$ to $\de_M^{\xi}$.
\end{proof}

\section{Discussion}

  We have shown that for an arbitrary finite poset
  the interval resolution global dimension ($\intgldim$) of its incidence algebra is finite (Proposition~\ref{prp:intgldim-fin})
  and provided a formula to compute $\intgldim$
  as the maximum of a finite number of terms (Proposition~\ref{prp:intdim-ta}).
  Beyond its applications in persistence, we hope that our results will be of independent interest in the representation theory of algebras.

  Removing some of the finiteness conditions may provide interesting
  (and useful) settings for further investigation.
  \begin{enumerate}
  \item What can be said about the resolution global dimension relative to other classes of indecomposable modules $\calX$ with an infinite number of elements?
    In the setting of finite posets, there are only a finite number of interval modules (up to isomorphism). 
    Our proof of the finiteness of $\intgldim$ (Proposition~\ref{prp:intgldim-fin}) relies heavily
    on this finiteness (from which follows that $X = G = \Ds_{I \in \bbI(\calP)} V_I$ is finitely generated)
    in order to apply Theorem~\ref{thm:ringel-qh} (cited from \cite[Theorem~in~\S 5]{ringel2010iyama}).
  \item What can be said in the setting of locally finite posets, or posets in general?
  We have obtained some conjectures about the value of the interval resolution global dimension
  of the $2$D commutative grids $\Gf{m,n}$ (Conjectures~\ref{cnj:ladder}~and~\ref{cnj:stabilize}).
  They also suggest that the interval resolution global dimension is not bounded
  as a function of $(m,n)$ if both $m$ and $n$ are allowed to vary. This suggests that a commutative grid infinite in both axes will have infinite $\intgldim$.
\end{enumerate}

In the
$2\times n$ commutative grid case, we provided
a new invariant (the modified compressed multiplicity) that is intimately related to an
alternating sum of the terms appearing in the interval resolution (Theorem~\ref{thm:comp-mult-resol}).
This is related to the result relating
the rank invariant to the projective resolutions in the rank-exact structure~\cite[Section~4]{botnan2021signed}.

Furthermore, the invariant obtained via M\"obius inversion of the modified compressed multiplicity
was shown to be a homological invariant (and thus a dim-hom invariant)
relative to the intervals (Corollary~\ref{cor:homologicalinvariant}).
Since the compressed multiplicity can be seen as a further modification of
the generalized rank invariant~\cite{kim2018generalized},
and since the generalized rank invariant was shown to not be a dim-hom invariant \cite[Corollary~7.10]{blanchette2021homological},
our results in the $2\times n$ case may point to potential modifications
to obtain a generalized rank invariant that is also a dim-hom invariant in the general case.

Extending our results for the $2\times n$ case to the general $m \times n$ case presents additional complications, especially in finding an appropriate morphism $\xi_I$ and generalizing Proposition~\ref{prp:L-int-dec}. At least, we were yet unable to find an immediate generalization. This may be a topic for further research.


\section*{Acknowledgements}
  On behalf of all authors, the corresponding author states that there is no conflict of interest.
  
In an initial version of this manuscript, we proved Proposition~\ref{prp:intdim-ta} directly. We thank an anonymous reviewer for suggesting the possibility of its generalization to Proposition~\ref{prp:I-res-gldim}.


\bibliography{refs}

\begin{thebibliography}{10}

\bibitem{asashiba2022interval}
Hideto Asashiba, Micka{\"e}l Buchet, Emerson~G Escolar, Ken Nakashima, and
  Michio Yoshiwaki.
\newblock On interval decomposability of 2d persistence modules.
\newblock {\em Computational Geometry}, page 101879, 2022.

\bibitem{asashiba2019approximation}
Hideto Asashiba, Emerson~G Escolar, Ken Nakashima, and Michio Yoshiwaki.
\newblock On approximation of $2$d persistence modules by
  interval-decomposables.
\newblock {\em arXiv preprint arXiv:1911.01637}, 2019.

\bibitem{assem2006elements}
Ibrahim Assem, Andrzej Skowro\'{n}ski, and Daniel Simson.
\newblock {\em Elements of the Representation Theory of Associative Algebras:
  Volume 1: Techniques of Representation Theory}, volume~65.
\newblock Cambridge University Press, 2006.

\bibitem{auslander1997representation}
Maurice Auslander, Idun Reiten, and Sverre~O. Smal\o.
\newblock {\em Representation theory of Artin algebras}, volume~36.
\newblock Cambridge university press, 1997.

\bibitem{bauer2020cotorsion}
Ulrich Bauer, Magnus~B Botnan, Steffen Oppermann, and Johan Steen.
\newblock Cotorsion torsion triples and the representation theory of filtered
  hierarchical clustering.
\newblock {\em Advances in Mathematics}, 369:107171, 2020.

\bibitem{blanchette2021homological}
Benjamin Blanchette, Thomas Br{\"u}stle, and Eric~J Hanson.
\newblock Homological approximations in persistence theory.
\newblock {\em arXiv preprint arXiv:2112.07632}, 2021.

\bibitem{botnan2020decomposition}
Magnus Botnan and William Crawley-Boevey.
\newblock Decomposition of persistence modules.
\newblock {\em Proceedings of the American Mathematical Society},
  148(11):4581--4596, 2020.

\bibitem{botnan2021signed}
Magnus~Bakke Botnan, Steffen Oppermann, and Steve Oudot.
\newblock Signed barcodes for multi-parameter persistence via rank
  decompositions and rank-exact resolutions.
\newblock {\em arXiv preprint arXiv:2107.06800}, 2021.

\bibitem{buchet_et_al:socg}
Micka{\"e}l Buchet and Emerson~G Escolar.
\newblock Realizations of indecomposable persistence modules of arbitrarily
  large dimension.
\newblock In {\em 34th International Symposium on Computational Geometry (SoCG
  2018)}. Schloss Dagstuhl-Leibniz-Zentrum fuer Informatik, 2018.

\bibitem{buchet2020every}
Micka{\"e}l Buchet and Emerson~G Escolar.
\newblock Every 1d persistence module is a restriction of some indecomposable
  2d persistence module.
\newblock {\em Journal of Applied and Computational Topology}, 4(3):387--424,
  2020.

\bibitem{carlsson2009theory}
Gunnar Carlsson and Afra Zomorodian.
\newblock The theory of multidimensional persistence.
\newblock {\em Discrete \& Computational Geometry}, 42(1):71--93, 2009.

\bibitem{crawley2015decomposition}
William Crawley-Boevey.
\newblock Decomposition of pointwise finite-dimensional persistence modules.
\newblock {\em Journal of Algebra and its Applications}, 14(05):1550066, 2015.

\bibitem{dey2021computing}
Tamal~K Dey, Woojin Kim, and Facundo M{\'e}moli.
\newblock Computing generalized rank invariant for 2-parameter persistence
  modules via zigzag persistence and its applications.
\newblock {\em arXiv preprint arXiv:2111.15058}, 2021.

\bibitem{dugas2007representation}
Alex~S Dugas.
\newblock Representation dimension as a relative homological invariant of
  stable equivalence.
\newblock {\em Algebras and representation theory}, 10(3):223--240, 2007.

\bibitem{edelsbrunner2008persistent}
Herbert Edelsbrunner and John Harer.
\newblock Persistent homology-a survey.
\newblock {\em Contemporary mathematics}, 453:257--282, 2008.

\bibitem{edelsbrunner2002topological}
Herbert Edelsbrunner, David Letscher, and Afra Zomorodian.
\newblock Topological persistence and simplification.
\newblock {\em Discrete Comput. Geom.}, 28(4):511--533, November 2002.

\bibitem{erdmann2004radical}
Karin Erdmann, Thorsten Holm, Osamu Iyama, and Jan Schr{\"o}er.
\newblock Radical embeddings and representation dimension.
\newblock {\em Advances in mathematics}, 185(1):159--177, 2004.

\bibitem{escolar2016persistence}
Emerson~G. Escolar and Yasuaki Hiraoka.
\newblock Persistence modules on commutative ladders of finite type.
\newblock {\em Discrete \& Computational Geometry}, 55(1):100--157, 2016.

\bibitem{gabriel1992algebra}
Peter Gabriel and Andrei~V. Roiter.
\newblock {\em Algebra VIII: Representations of Finite-dimensional Algebras},
  volume~73.
\newblock Springer Verlag, 1992.

\bibitem{iyama2003finiteness}
Osamu Iyama.
\newblock Finiteness of representation dimension.
\newblock {\em Proceedings of the american mathematical society},
  131(4):1011--1014, 2003.

\bibitem{kim2018generalized}
Woojin Kim and Facundo M{\'e}moli.
\newblock Generalized persistence diagrams for persistence modules over posets.
\newblock {\em Journal of Applied and Computational Topology}, 5(4):533--581,
  2021.

\bibitem{leszczynski1994representation}
Zbigniew Leszczy{\'n}ski.
\newblock On the representation type of tensor product algebras.
\newblock {\em Fundamenta Mathematicae}, 144(2):143--161, 1994.

\bibitem{leszczynski2000tame}
Zbigniew Leszczynski and Andrzej Skowronski.
\newblock Tame triangular matrix algebras.
\newblock In {\em Colloq. Math}, volume~86, pages 259--303, 2000.

\bibitem{ringel2010iyama}
Claus~Michael Ringel.
\newblock Iyama’s finiteness theorem via strongly quasi-hereditary algebras.
\newblock {\em Journal of Pure and Applied Algebra}, 214(9):1687--1692, 2010.

\bibitem{rota1964foundations}
Gian-Carlo Rota.
\newblock On the foundations of combinatorial theory i. theory of m{\"o}bius
  functions.
\newblock {\em Probability theory and related fields}, 2(4):340--368, 1964.

\bibitem{GAP4}
{The GAP~Group}.
\newblock {GAP -- Groups, Algorithms, and Programming, Version 4.11.1}.
\newblock \url{https://www.gap-system.org}, 2021.

\bibitem{qpa}
{The QPA-team}.
\newblock {QPA - Quivers, path algebras and representations - a GAP package,
  Version 1.31}.
\newblock \url{https://folk.ntnu.no/oyvinso/QPA/}, 2020.

\end{thebibliography}

\end{document}